\documentclass[12pt]{article}

\usepackage{amsmath,amssymb,enumerate,bbm}
\newtheorem{theorem}{Theorem}[section]

\newtheorem{definition}[theorem]{Definition}

\numberwithin{equation}{section}

 \makeatletter
 
 \newcommand{\Rmnum}[1]{\expandafter\@slowromancap\romannumeral #1@}
 \makeatother

\begin{document}

\title{A construction of $3$-e.c. graphs using quadrances}\author{Le Anh
Vinh\\
Mathematics Department\\
Harvard University\\
Cambridge, MA 02138, US\\
vinh@math.harvard.edu}\maketitle

\begin{abstract}
  A graph is $n$-e.c. ($n$-existentially closed) if for every pair of subsets
  $A, B$ of vertex set $V$ of the graph such that $A \cap B = \emptyset$ and
  $|A| + |B| = n$, there is a vertex $z$ not in $A \cup B$ joined to each
  vertex of $A$ and no vertex of $B$. Few explicit families of $n$-e.c. are
  known for $n > 2$. In this short note, we give a new construction of
  $3$-e.c. graphs using the notion of quadrance in the finite Euclidean space
  $\mathbbm{Z}_p^d$.
\end{abstract}

\section{Introduction}

For a positive integer $n$, a graph is $n$-\textit{existentially closed} or
$n$-e.c. if we can extend all $n$-subsets of vertices in all possible ways.
Precisely, if for every pair of subsets $A, B$ of vertex set $V$ of the graph
such that $A \cap B = \emptyset$ and $|A| + |B| = n$, there is a vertex $z$
not in $A \cup B$ joined to each vertex of $A$ and no vertex of $B$. From the
results of Erd\"os and R\'enyi \cite{er}, almost all finite graphs are $n$-e.c. Despite this result, until recently, only few explicit examples of $n$-e.c. graphs are
known for $n > 2$ (see \cite{anthony} for a comprehensive survey on the constructions of $n$-e.c. graphs). In this short note, we give a new construction of $3$-e.c. graphs
using the notion of quadrance in the finite Euclidean space $\mathbbm{Z}_p^d$.

Suppose that $p$ be an odd prime, and that $\mathbbm{Z}_p =\{0, \ldots, p -
1\}$ be the prime field with $p$ elements. We will construct a $3$-e.c. graph
with the vertex set $\mathbbm{Z}_p^d$ for some large $d$. The following
definition of quadrance is taken from \cite{norman}. 

\begin{definition}
  The quadrance between the points $X = (x_1, \ldots, x_d)$ and $Y
  (y_1, \ldots, y_d)$ in $\mathbbm{Z}_p^d$ is the number
  \[ Q (X, Y) := (x_1 - y_1)^2 + \ldots + (x_d - y_d)^2 \in \mathbbm{Z}_p
     . \]
\end{definition}

Let $V_1 =\{0, 1, 2, \ldots, (p - 1) / 2\}$. We define the graph $G_{p, d}$ as
follows. The vertices of the graph $G_{p, d}$ are the points of
$\mathbbm{Z}_p^d$. There is an edge between two vertices $X$ and $Y$ if and
only if $Q (X, Y) \in V_1$. We claim that $G_{p, d}$ is $3$-e.c. for $p
\geqslant 7$ and $d \geqslant 5$.

\begin{theorem} \label{main}
  Suppose that $p \geqslant 7$ be an odd prime and $d \geqslant 5$ be an
  integer. Then the graph $G_{p, d}$ is $3$-e.c.
\end{theorem}

Note that these quadrance graphs are just Cayley graphs of $\mathbbm{Z}_p^d$.

\section{The $3$-e.c. property of the graph $G_{p,d}$}

  We now give a proof of Theorem \ref{main}.
  Let $V_2 =\{(p+1)/2 \ldots, p - 1\}=\mathbbm{Z}_p \backslash V_1$. It suffices
  to show that for any three distinct points $A = (a_1, \ldots, a_d)$, $B =
  (b_1, \ldots, b_d)$, $C = (c_1, \ldots, c_d)$ in $\mathbbm{Z}_p^d$ and $i,
  j, k \in \{1, 2\}$, there is a point $X = (x_1, \ldots, x_d) \in
  \mathbbm{Z}_p^d$, $X \neq A, B, C$ such that $Q (X, A) \in V_i$, $Q (X, B)
  \in V_j$ and $Q (X, C) \in V_k$. Therefore, we only need to show that there
  exist $u \in V_i, v \in V_j$, and $w \in V_k$ such that the following system
  has at least four solutions (in this case, one of these solutions is different from $A$, $B$, and $C$),
  \begin{eqnarray}
    (x_1 - a_1)^2 + \ldots + (x_d - a_d)^2 & = & u\\
    (x_1 - b_1)^2 + \ldots + (x_d - b_d)^2 & = & v\label{e2}\\
    (x_1 - c_1)^2 + \ldots + (x_d - c_d)^2 & = & w.\label{e3}
  \end{eqnarray}
  For any $ X = (x_1,\ldots,x_d) \in \mathbbm{Z}_p^d$, define
  \[ \| X\| = x_1^2 + \ldots + x_d^2.\]
  By eliminating $x_i^2$'s from (\ref{e2}) and (\ref{e3}), we get an equivalent system of equations
  \begin{eqnarray}
    Q (X, A) & = & u\label{e4}\\
    \left\langle X, B - A \right\rangle & = & (u - v +\|B\|-\|A\|) / 2 \label{e5}\\
    \left\langle X, C - A \right\rangle & = & (u - w +\|C\|-\|A\|) / 2. \label{e6}
  \end{eqnarray}
  We first show that the system of two equations (\ref{e5}) and (\ref{e6}) has a solution $X_0$ for some choices of $u \in V_i$, $v \in V_j$, and $w \in V_k$. We consider two cases.
  
  Case 1. Suppose that $B - A$ and $C - A$ are linearly independent. For
  any $u \in V_i$, $v \in V_j$, and $w \in V_k$, it is clear that there is a solution $X_0$ to the system of two equations (\ref{e5}) and (\ref{e6}).
  
  Case 2. Suppose that  $B - A$ and $C - A$ are linearly dependent. Since $C - A \neq B - A \neq 0$, $C - A = t (B - A)$ for some $t \neq 0, 1$. The two equations (\ref{e5}) and (\ref{e6}) have a common solution if we can choose $u \in V_i$, $v \in V_j$, and $w \in V_k$ such that
  \[ u - w +\|C\|-\|A\|= t (u - v +\|B\|-\|A\|), \]
  or equivalently,
  \[ w = t v + a, \]
  where $a =\|C\|+ (t - 1)\|A\|- t\|B\|- (t - 1) u$. In other words, we need
  to show that $\{t v : v \in V_j \} \cap \{w - a : w \in V_k \} \neq
  \emptyset$. We have two subcases.
  
\begin{itemize}
	\item Suppose that $t \neq 0, \pm 1$. We label $\mathbbm{Z}_p$ around the
  circle. The set $\{w - a : w \in V_k \}$ is a block of $(p \pm 1) / 2$
  consecutive points. Going $|V_k | = (p \pm 1) / 2$ steps of length $2 < |t|
  \leqslant (p - 1) / 2$ around the circle, we cannot avoid any block of $(p
  \pm 1) / 2$ consecutive points. Hence, for any fixed $u \in V_i$, we can  choose $v \in V_j$ and $w \in V_k$ such that $w = t v + a$. 
  \item Suppose that $t = - 1$. The set $\{w + v : w \in V_k, v \in V_j \}$ contains at least $p - 2$ elements. Since $|A_i | \geqslant (p - 1) / 2 \geqslant 3$, we can choose
  $u$ such that $a \in \text{$\{w + v : w \in V_k, v \in V_j \}$}$.
\end{itemize}
  Therefore, we always can choose $u \in V_i$, $v \in V_j$, and $w \in V_k$ such that the two equations (\ref{e5}) and (\ref{e6}) have a common solution $X_0$.
  
  Take a basis of solutions of the system
  \begin{eqnarray*}
    \left\langle X, B - A \right\rangle & = & 0\\
    \left\langle X, C - A \right\rangle & = & 0,
  \end{eqnarray*}
  and the solution $X_0$. Substitute them into (\ref{e4}), we get a single
  quadratic equation of $d - 2$ variables. Since $d - 2 \geqslant 3$, this
  quadratic equation has at least $p$ $(\geq 4)$ solutions. Theorem \ref{main} follows immediately.
  
\section{Remarks and Further Questions}

Note that the construction is well defined over $\mathbbm{Z}_m$ for any $m \in \mathbbm{N}$ and it gives $3$-e.c. graphs as well. The proof goes without any essential changes when $p$ is not a prime. 

Moreover, the proof of Theorem \ref{main} only works for $d\geq 5$. It is plausible to conjecture that the graphs are $3$-e.c. for $d \geq 2$. Another interesting question is to consider other constructions with difference choices of $V_1 \subset \mathbbm{Z}_p$. When $d=2$, let $V =\{a^2 : a\in \mathbbm{Z}_p^* \}$. We define the graph $G_{V,p}$ as
follows. The vertices of the graph $G_{V, p}$ are the points of
$\mathbbm{Z}_p^2$. There is an edge between two vertices $X$ and $Y$ if and
only if $Q (X, Y) \in V$. We know that $G_{V,p}$ is isomorphic to the Paley graph $P_p$ (see, for example, \cite{vinh}). It is well known that $P_p$ is $n$-e.c for any $n$ given that $p$ is sufficiently large, so is $G_{V,p}$. We, however, have not known any results for the remaining cases. 

\section*{Acknowledgment}

The author would like to thank Prof. Igor Shparlinski for helpful advice on the proof of Theorem \ref{main}. He also wants to thank the referee for constructive comments and suggestions.

\end{document}